\documentclass[12pt]{report}
\usepackage{amssymb,amsmath,makeidx,verbatim}
\newcommand{\R}{\mathbb{R}}

\newcommand{\C}{\mathbb{C}}

\newcommand{\be}{\begin{enumerate}}
\newcommand{\ee}{\end{enumerate}}
\newcommand{\bq}{\begin{eqnarray*}}
\newcommand{\eq}{\end{eqnarray*}}

\begin{document}
\newcommand{\disp}{\displaystyle}
\thispagestyle{empty}
\begin{center}
\textsc{Abstract Peter-Weyl theory for semicomplete orthonormal sets.\\}
\ \\
\textsc{Olufemi O. Oyadare}\\
\ \\
Department of Mathematics,\\
Obafemi Awolowo University,\\
Ile-Ife, $220005,$ NIGERIA.\\
\text{E-mail: \textit{femi\_oya@yahoo.com}}\\
\end{center}
\begin{quote}
{\bf Abstract.} {\it The central concept in the harmonic analysis of a compact group is the completeness of Peter-Weyl orthonormal basis as constructed from the matrix coefficients of a maximal set of irreducible unitary representations of the group, leading ultimately to the direct sum decomposition of its $L^{2}-$ space. A Peter-Weyl theory for a semicomplete orthonormal set is also possible and is here developed in this paper for compact groups. Existence of semicomplete orthonormal sets on a compact group is proved by an explicit construction of the standard Riemann-Lebesgue semicomplete orthonormal set. This approach gives an insight into the role played by the $L^{2}-$ space of a compact group, which is discovered to be just an example (indeed the largest example for every semicomplete orthonormal set) of what is called a prime-Parseval subspace, which we proved to be dense in the usual $L^{2}-$ space, serves as the natural domain of the Fourier transform and breaks up into a direct-sum decomposition. This paper essentially gives the harmonic analysis of the prime-Parseval subsapce of a compact group corresponding to any semicomplete orthonormal set, with an introduction to what is expected for all connected semisimple Lie groups through the notion of a $K-$semicomplete orthonormal set.}
\end{quote}
\ \\
\ \\
\ \\
\ \\
$\overline{2010\; \textmd{Mathematics}}$ Subject Classification: $43A85, \;\; 22E30, \;\; 22E46$
\ \\
\ \\
\ \\
{\bf \S 1. Introduction.}

Harmonic analysis on a compact group is mainly a direct consequence of the famous \textit{Peter-Weyl theory} which gives a consistent method, via the computation of the matrix coefficients of its \textit{irreducible unitary representations,} of deriving a \textit{complete orthonormal set} which is immediately responsible for the direct-sum decomposition of its $L^{2}-$ space and \textit{regular representation.} Even though such a complete orthonormal set is non-existence for non-compact topological groups and hence the harmonic analysis on \textit{non-compact topological groups,} as we know for \textit{connected nilpotent} and \textit{semisimple Lie groups,} has had to be developed through other means notably via the differential equations satisfied by the (\textit{spherical}) functions derived as matrix coefficients of irreducible unitary representations constructed from \textit{parabolic} and \textit{cohomological inductions} and the completeness afforded by the \textit{Plancherel theorem} (which in the final analysis still depends on the availability and properties of the \textit{discrete series} (known to be the irreducible unitary representations corresponding to some complete orthonormal set) of some distinguished compact subgroups), it still found to be appropriate (and to have a sense of finality) to have some forms of \textit{Peter-Weyl} results on such \textit{non-compact topological groups.}

It is however possible to get at the decomposition of the regular representation of a compact group $G$ (for a start) via the indirect use of the notion of a \textit{semicomplete orthonormal set} on such a group, leading to the consideration of a distinguished subspace of $L^{2}(G)$ which is established to be \textit{topologically dense.} The study in this paper opens up this field of research by a detailed look at the \textit{compact case.} The paper is arranged as follows.

$\S 2.$ contains a quick review of the well-known notion of a complete orthonormal set on a compact group, giving the detailed of the aforementioned consistent way of constructing such a set through \textit{Peter-Weyl theorem} which then leads to the direct-sum decomposition of its $L^{2}-$space. The concept of a \textit{semicomplete orthonormal set} on a compact group $G$ is introduced in $\S 3.$ with constructible examples (prominent among which is the \textit{Riemann-Lebesgue} orthonormal set), where we derived and used the properties of the \textit{Fourier} and \textit{prime-Parseval subspaces} of $L^{2}(G).$ Chief among these properties is the topological denseness of every \textit{prime-Parseval subspace} in $L^{2}(G).$ This takes us to the \textit{Fourier transform} of the \textit{prime-Parseval subspace} and its direct-sum decomposition into \textit{invariant subspaces.} The last section gives an introductory extension of the results of $\S 3.$ on compact groups to connected semisimple Lie groups with finite center.
\ \\
\ \\
{\bf \S 2. Fourier and Parseval subsapces for complete orthonormal set.}

A mutually orthonormal family $\{\chi_{\alpha}\}_{\alpha\in A}$ in a Hilbert space, $(H,\langle\cdot,\cdot\rangle)$ is said to be \textit{complete} (in $H$) if $x\in H$ is such that $\langle x,\chi_{\alpha}\rangle=0$ (for every $\alpha\in A$) implies $x=0.$ This means that a family $\{\chi_{\alpha}\}_{\alpha\in A}$ of mutually orthonormal members of $H$ is complete whenever it can be shown that the zero element of $H$ is the \textit{only} non-member of the family that is mutually orthonormal to all members of the said family. Two other equivalent methods of confirming the completeness of the family $\{\chi_{\alpha}\}_{\alpha\in A}$ are as follows.

\textbf{2.1 Lemma.} ([5.], p.\;3) \textit{Let $\{\chi_{\alpha}\}_{\alpha\in A}$ denote a mutually orthonormal family in a Hilbert space $(H,\langle\cdot,\cdot\rangle).$ The following are equivalent:\\
(a) Every $x\in H$ can be expressed as $x=\sum_{\alpha\in A}\langle x,\chi_{\alpha}\rangle\chi_{\alpha}.$\\
(b) Every $x\in H$ satisfies $\parallel x\parallel^{2}=\sum_{\alpha\in A}\mid\langle x,\chi_{\alpha}\rangle\mid^{2}.$\\
(c) $\{\chi_{\alpha}\}_{\alpha\in A}$ is \textit{complete} in $H.\;\Box$}

The informed reader would observe that $(a)$ of $(2.1)$ is a \textit{Fourier series} expansion of $x$ while $(b)$
of $(2.1)$ is its \textit{Parseval equality,} both with respect to $\{\chi_{\alpha}\}_{\alpha\in A}.$ The import of this equivalence (in the light of $(a)$ of $(2.1)$ (respectively, $(b)$ of $(2.1)$)) is that every $x\in H$ has a Fourier series expansion in terms of any known complete orthonormal family in $H.$ We could then say that the subset $H(\chi_{\alpha})$ of $H$ given as $$\{x\in H: x=\sum_{\alpha\in A}\langle x,\chi_{\alpha}\rangle\chi_{\alpha},\;\mbox{for some orthonormal family}\; \{\chi_{\alpha}\}_{\alpha\in A}\;\mbox{in}\;H\}$$ (equivalently, the subset $H_{\mathfrak{P}}(\chi_{\alpha})$ of $H$ given also as $$\{x\in H: \parallel x\parallel^{2}=\sum_{\alpha\in A}\mid\langle x,\chi_{\alpha}\rangle\mid^{2},\;\mbox{for some orthonormal family}\; \{\chi_{\alpha}\}_{\alpha\in A}\;\mbox{in}\;H\})$$ is exactly $H$ if, and only if, $\{\chi_{\alpha}\}_{\alpha\in A}$ is complete. Indeed another version of the equivalence of Lemma $2.1,$ whose formulation serves as our point of departure, is given as follows.

\textbf{2.2 Lemma.} \textit{Let $\{\chi_{\alpha}\}_{\alpha\in A}$ denote a mutually orthonormal family in a Hilbert space $(H,\langle\cdot,\cdot\rangle).$ The following are equivalent:\\
(a) $H(\chi_{\alpha})=H$\\
(b) $H_{\mathfrak{P}}(\chi_{\alpha})=H$\\
(c) $\{\chi_{\alpha}\}_{\alpha\in A}$ is \textit{complete} in $H.\;\Box$}

\textbf{2.3 Remarks.} It may be safely conjectured that the \textit{Fourier subspace} $H(\chi_{\alpha})$ as well as the \textit{Parseval subspace} $H_{\mathfrak{P}}(\chi_{\alpha})$ (of a Hilbert space $H$) with respect to a complete mutually orthonormal family will always be equal to $H.$ It will be a delight to study the disparity between the \textit{Fourier subspace} $H(\chi_{\alpha})$ as well as the \textit{Parseval subspace} $H_{\mathfrak{P}}(\chi_{\alpha})$ (of $H$ with respect to the mutually orthonormal family $\{\chi_{\alpha}\}_{\alpha\in A}$) and their inclusions in $H,$ when the family $\{\chi_{\alpha}\}_{\alpha\in A}$ is not complete.

For example, if the family $\{\chi_{\alpha}\}_{\alpha\in A}$ of mutually orthonormal members in $H$ is such that \textit{$\langle x,\chi_{\alpha}\rangle=0$ (for every $\alpha\in A$) does not necessarily imply whether $x=0$ or $x\neq0,$} it possible to then have that $$0\leq \parallel x\parallel^{2}=\sum_{\alpha\in A}\mid\langle x,\chi_{\alpha}\rangle\mid^{2}=0,$$ showing in this case (for the family $\{\chi_{\alpha}\}_{\alpha\in A}$ in which $\langle x,\chi_{\alpha}\rangle=0$ (for every $\alpha\in A$) does not necessarily imply whether $x=0$ or $x\neq0$) that we now have $H_{\mathfrak{P}}(\chi_{\alpha})=\{0\}$ ($=H(\chi_{\alpha})\neq H,$ showing that both subspaces are too small and far from being equal to $H$). This shows at a glance the importance of completeness of the family $\{\chi_{\alpha}\}_{\alpha\in A}$ in the consideration of the \textit{Parseval equality}, for the non-triviality of these two subspaces $H(\chi_{\alpha})$ and $H_{\mathfrak{P}}(\chi_{\alpha})$ and for the sustenance of the relationship of equality (of Lemma $2.2$) between $H(\chi_{\alpha})$ and $H_{\mathfrak{P}}(\chi_{\alpha}).\;\Box$

However, and as it shall be shown in the next section, these two subspaces, $H(\chi_{\alpha})$ and $H_{\mathfrak{P}}(\chi_{\alpha})$ may be considered for an \textit{appropriately chosen} not-necessarily complete orthonormal family $\{\chi_{\alpha}\}_{\alpha\in A}$ and with which they would still be found not to be too small in sizes (in comparison with $H$). This choice of a not-necessarily complete orthonormal family $\{\chi_{\alpha}\}_{\alpha\in A}$ would equally help and be appropriate in order that both $H(\chi_{\alpha})$ and $H_{\mathfrak{P}}(\chi_{\alpha})$ be \textit{lifted} to all of $H.$ All this in a moment.

A well-known method of computing complete orthonormal family of functions is via the matrix coefficients of irreducible unitary representations of a compact groups $G$ which is then used to decompose $L^{2}(G)$ into invariant subspaces, leading to the decomposition of the right regular representation on $G$ (which sadly, does not generalize to \textit{non-compact topological groups}). Here is the technique.

Denote the \textit{dual} of a compact group $G$ by $\widehat{G},$ consisting of all its equivalence classes of irreducible unitary representations. For $\lambda\in\widehat{G}$ denote by $u_{ij}^{\lambda}$ the corresponding matrix coefficient representative of the class $\lambda$ whose degree is also denoted by $d(\lambda).$ Then the set $$\{\sqrt{d(\lambda)}u_{ij}^{\lambda}:\;\lambda\in\widehat{G},\;1\leq i,j\leq d(\lambda)\}$$ consists of a maximal set of complete orthonormal family of functions in $L^{2}(G)$ and (hence) every $f\in L^{2}(G)$ can be expanded as $$f=\sum_{\lambda\in\widehat{G}}d(\lambda)\sum_{i,j}^{d(\lambda)}\langle f,u_{ij}^{\lambda}\rangle u_{ij}^{\lambda}$$ (with convergence in the norm of $L^{2}(G)$) whose \textit{Fourier transform} $$\widehat{f}:\widehat{G}\rightarrow M_{d(\lambda)}(\mathbb{C}):\lambda\mapsto\widehat{f}(\lambda)=(\widehat{f}(\lambda)_{ij})_{i,j=1}^{d(\lambda)}$$ is given as $\widehat{f}(\lambda)_{ij}:=\langle f,u_{ij}^{\lambda}\rangle$ (where $M_{d(\lambda)}(\mathbb{C})$ denotes the algebra of matrices with entries in $\mathbb{C}$ and degree $d(\lambda)$). It then follows that for any compact group $G,$ the \textit{Fourier subspace} $L^{2}(G)(\sqrt{d(\lambda)}u_{ij}^{\lambda})$ of $L^{2}(G)$ is given as $$L^{2}(G)(\sqrt{d(\lambda)}u_{ij}^{\lambda}):=\{f\in L^{2}(G):\;f=\sum_{\lambda\in\widehat{G}}d(\lambda)\sum_{i,j}^{d(\lambda)}\langle f,u_{ij}^{\lambda}\rangle u_{ij}^{\lambda}\}=L^{2}(G)$$ (=$L^{2}(G)_{\mathfrak{P}}(\sqrt{d(\lambda)}u_{ij}^{\lambda}),$ the \textit{Parseval subspace} of $L^{2}(G)$), with respect to the family $\{\sqrt{d(\lambda)}u_{ij}^{\lambda}:\;\lambda\in\widehat{G},\;1\leq i,j\leq d(\lambda)\}.$ We then have the abstract direct-sum decomposition of $L^{2}(G)$ given as $$L^{2}(G)=\bigoplus_{\lambda\in\widehat{G}}\bigoplus_{i=1}^{d(\lambda)}H_{i}^{\lambda},$$ where $H_{i}^{\lambda}:=\sum_{j=1}^{d(\lambda)}\mathbb{C}u_{ij}^{\lambda}.$ This is the content of \textit{Peter-Weyl Theorem,} $[5.],$ and we shall refer to the set $$\{\sqrt{d(\lambda)}u_{ij}^{\lambda}:\;\lambda\in\widehat{G},\;1\leq i,j\leq d(\lambda)\}$$ as the \textit{standard Peter-Weyl orthonormal set} on $G.$

The inability of being able to get an orthonormal family in $L^{2}(G)$ for a non-compact topological group $G$ in the above tradition of Peter-Weyl is the first stumbling block to harmonic analysis on such groups, which has been considerably understood and completely developed via a rigorous treatment of the rich structure of differential equations satisfied by matrix-coefficients of members of each of the classes in $\widehat{G},\;[2].$ This paper presents a constructive method of getting a not-necessarily complete orthonormal set which is \textit{close enough} to being a complete orthonormal family in an arbitrary Hilbert space $(H,\langle\cdot,\cdot\rangle)$ and/or in $L^{2}(G),$ for a compact group (and introduced the same technique for a semisimple Lie group) offering a more general Fourier series expansion of each member of an appropriate subspace of $H$ and/or $L^{2}(G).$

Starting with a compact group (before extending the notion to all connected semisimple Lie groups, with finite center, via its \textit{Iwasawa decomposition}) we would however not approach harmonic analysis on the groups via the completeness (and consequent denseness) of the \textit{standard Peter-Weyl orthonormal set,} but via a denseness in the $L^{2}-$space which would be found to be possible from an \textit{almost complete} orthonormal set.
\ \\
\ \\
{\bf \S 3. Semicomplete orthonormal set in a compact group.}

The existence of different special functions and polynomials of mathematical physics, which have been established to be orthonormal in various semisimple Lie groups (compact and non-compact types), is well-known. However the absence of completeness of these orthornormal families (under the structure of their individual corresponding groups) is the first stumbling block to a direct \textit{Peter-Weyl harmonic analysis} of them. In this section we shall define and study the concept of a \textit{semicomplete} orthonormal family in a compact group in order to extend this concept to the harmonic analysis of \textit{all} semisimple Lie groups in the next section.

\textbf{3.1 Definition.} (\textit{Semicomplete orthonormal family}) \textit{Let $G$ denote a compact group and let the members of the non-empty set $A$ be ordered such that $A=\{\alpha_{i}^{j}\}_{i,j}.$ An orthonormal family $\{\chi_{\alpha_{i}^{j}}\}_{\alpha_{i}^{j}\in A}$ in $L^{2}(G)$ is said to be semicomplete if given $\epsilon>0$ there exist some non-zero scalars $$\gamma_{1},\cdots,\gamma_{k},\cdots,\beta_{11},\cdots,\beta_{ij},\cdots\in\mathbb{C}$$ and $n\in\mathbb{N}$ such that $$\parallel\sum_{\lambda\in\widehat{G}}d(\lambda)\sum_{i,j=1}^{d(\lambda)}\langle f,u_{ij}^{\lambda}\rangle u_{ij}^{\lambda}-\sum_{j=1}^{n}\gamma_{j}\sum_{i=1}^{n}\beta_{ij}\sum_{\alpha_{i}^{j}\in A}\langle f,\chi_{\alpha_{i}^{j}}\rangle\chi_{\alpha_{i}^{j}}\parallel_{2}<\epsilon$$ for every $f\in L^{2}(G).\;\Box$}

The quantity $$\sum_{\lambda\in\widehat{G}}d(\lambda)\sum_{i,j=1}^{d(\lambda)}\langle f,u_{ij}^{\lambda}\rangle u_{ij}^{\lambda}$$ in Definition $3.1$ above may be replaced with $f$ (due to the \textit{Peter-Weyl Theorem}), so that the other quantity $$\sum_{j=1}^{n}\gamma_{j}\sum_{i=1}^{n}\beta_{ij}\sum_{\alpha_{i}^{j}\in A}\langle f,\chi_{\alpha_{i}^{j}}\rangle\chi_{\alpha_{i}^{j}}$$ (in the same Definition above) should be seen as the \textit{total contribution of $\{\chi_{\alpha_{i}^{j}}\}_{\alpha_{i}^{j}\in A}$ in $L^{2}(G)$ in its bid to attain $f.$} Thus the informed reader would see that the inequality in Definition $3.1$ above simply gives a measure of how close to the completeness (of $\{\sqrt{d(\lambda)}u_{ij}^{\lambda}\}$) is the orthonormal set $\{\chi_{\alpha_{i}^{j}}\}_{\alpha_{i}^{j}\in A}.$

The \textit{standard Peter-Weyl orthonormal basis} $\{\sqrt{d(\lambda)}u_{ij}^{\lambda}\}$ used in the above Definition $3.1$ may be replaced by any other known complete orthonormal set $\{v_{\mu}\}_{\mu\in B}$ in $L^{2}(G)$ while the concept of a semicomplete orthonormal set (for $\{\chi_{\alpha_{i}^{j}}\}_{\alpha_{i}^{j}\in A}$) could also be defined for an arbitrary Hilbert space, $H,$ so as to have what may be generally called \textit{a semicomplete orthonormal set in $H$ with respect to (the complete orthonormal set) $\{v_{\mu}\}_{\mu\in B}$ in $H.$} If in this general case the set $\{v_{\mu}\}_{\mu\in B}$ in $H$ is also not necessarily complete, we may arrive at the notion of \textit{a relative semicomplete orthonormal set} for $\{\chi_{\alpha_{i}^{j}}\}_{\alpha_{i}^{j}\in A}$ in $H$ with respect to $\{v_{\mu}\}_{\mu\in B}$ in $H.$ Thus Definition $3.1$ may therefore be seen as giving \textit{semicompleteness of $\{\chi_{\alpha_{i}^{j}}\}_{\alpha_{i}^{j}\in A}$ in $L^{2}(G)$ with respect to the standard Peter-Weyl orthonormal basis $\{\sqrt{d(\lambda)}u_{ij}^{\lambda}\}.$}

It is clear that every complete orthonormal set in $L^{2}(G)$ (or in any Hilbert space, $H$) is automatically semicomplete; simply choose $A=\widehat{G},\;\gamma_{j}=\beta_{ij}=1,$ but not conversely. An inductive method of immediately constructing a semicomplete orthonormal set in a compact group is by a method of \textit{selective omission} of some number of members in any known complete (or of the \textit{standard Peter-Weyl}) orthonormal set with a \textit{controlled bound.} The control of the bound in the method of \textit{selective omission} would be achieved using the \textit{Riemann-Lebesgue Lemma.}

This method, as contained in the following, equally gives an \textit{existence} argument for the concept of a semicomplete orthonormal set in a compact group.

\textbf{3.2 Lemma.} (\textit{Existence of a semicomplete orthonormal set: the standard Riemann-Lebesgue orthonormal set}) \textit{Let $G$ denote a compact group. Then there exist $\lambda_{0}\in \widehat{G}$ for which $$\mid\langle f,u_{km}^{\lambda}\rangle\mid<\frac{\epsilon}{d(\lambda_{0})},$$ for every $f\in L^{2}(G),\;\mid\lambda\mid\geq\mid\lambda_{0}\mid$ and $1\leq k,m\leq d(\lambda_{0}).$ Moreover, $$\{\sqrt{d(\lambda)}u_{ij}^{\lambda}:\;\lambda\in\widehat{G}\setminus\{\lambda_{0}\},\;1\leq i,j\leq d(\lambda)\}$$ is a semicomplete orthonormal set on the compact group.}

\textbf{Proof.} Since the dual group $\widehat{G}$ of a compact group $G$ is discrete, so that $$\lim_{\mid\lambda\mid\rightarrow\infty}\langle f,u_{ij}^{\lambda}\rangle=\lim_{\mid\lambda\mid\rightarrow\infty}\widehat{f(\lambda)}_{ij}=0\;\;(\mbox{by the \textit{Riemann-Lebesgue Lemma}}),$$ it follows that there are (infinitely) many possible $\lambda\in \widehat{G}$ (choose such one $\lambda_{0}$) with $\mid\lambda\mid\geq\mid\lambda_{0}\mid$ for which $\mid\langle f,u_{km}^{\lambda}\rangle\mid=\mid\langle f,u_{km}^{\lambda}\rangle-0\mid<\frac{\epsilon}{d(\lambda_{0})},$ for every $f \in L^{2}(G)$ and  $1\leq k,m\leq d(\lambda_{0}),$ as required. Hence,
$$\parallel\sum_{\lambda\in\widehat{G}}d(\lambda)\sum_{i,j=1}^{d(\lambda)}\langle f,u_{ij}^{\lambda}\rangle u_{ij}^{\lambda}-\sum_{\lambda\in\widehat{G}\setminus\{\lambda_{0}\}}d(\lambda)\sum_{i,j=1}^{d(\lambda)}\langle f,u_{ij}^{\lambda}\rangle u_{ij}^{\lambda}\parallel_{2}=\parallel d(\lambda_{0})\sum_{i,j=1}^{d(\lambda_{o})}\langle f,u_{ij}^{\lambda_{0}}\rangle u_{ij}^{\lambda_{0}}\parallel_{2}$$ $\leq d(\lambda_{0})\sum_{i,j=1}^{d(\lambda_{0})}\mid\langle f,u_{ij}^{\lambda_{0}}\rangle\mid<\epsilon,$ for every $f\in L^{2}(G).\;\Box$

The technique of Lemma $3.2$ may be extended as follows. Generally, choose (as assured by the \textit{Riemann-Lebesgue Lemma}) $\lambda_{0}^{(1)},\lambda_{0}^{(2)},\cdots\in\widehat{G}$ for which $$\sum_{k=1}^{\infty}\mid\langle f,u_{ij}^{\lambda}\rangle\mid<\frac{\epsilon}{\sum_{k=1}^{\infty}d(\lambda_{0}^{(k)})}$$ where $\mid\lambda\mid\geq\max\{\mid\lambda_{0}^{(1)}\mid,\mid\lambda_{0}^{(2)}\mid,\;\cdots\}$ and $f \in L^{2}(G).$ Then, with proof essentially the same as in Lemma $3.2,$ the set $$\{\sqrt{d(\lambda)}u_{ij}^{\lambda}:\;\lambda\in\widehat{G}\setminus\{\lambda_{0}^{(1)},\lambda_{0}^{(2)},\cdots\},\;1\leq i,j\leq d(\lambda)\}$$ is a semicomplete orthonormal set on the compact group, $G.$ We shall henceforth refer to the semicomplete orthonormal set $$\{\sqrt{d(\lambda)}u_{ij}^{\lambda}:\;\lambda\in\widehat{G}\setminus\{\lambda_{0}^{(1)},\lambda_{0}^{(2)},\cdots\},\;1\leq i,j\leq d(\lambda)\}$$ as the \textit{standard Riemann-Lebesgue (semicomplete) orthonormal set} (being in correspondence with the \textit{standard Peter-Weyl (complete) orthonormal set,} $\{\sqrt{d(\lambda)}u_{ij}^{\lambda}\}.$)

Other \textit{non-standard} examples of Definition $3.1$ may be deduced from the numerous \textit{special functions} of mathematical physics where their corresponding non-zero scalars $\gamma_{j}$ and $\beta_{ij}$ in Definition $3.1$ could be calculated from.

\textbf{3.3 Remarks.} In contrast to the zero-subspace $H_{\mathfrak{P}}(\chi_{\alpha})$ of Remarks $2.3$ we may, in the context of a semicomplete orthonormal set $\{\chi_{\alpha}\}_{\alpha\in A}$ in a Hilbert space $(H,\langle\cdot,\cdot\rangle),$ consider the subspace $$H_{\mathfrak{P}}'(\chi_{\alpha}):=\{x\in H:\;\langle x,\chi_{\alpha}\rangle=0,\;\mbox{(for every $\alpha\in A$) implies}\;x=0\},$$ for some orthonormal set $\{\chi_{\alpha}\}_{\alpha\in A}$ in $H.$ It is clear (from Lemma $2.2$) that $H_{\mathfrak{P}}'(\chi_{\alpha})=H$ (hence equal to $H(\chi_{\alpha})$ and $H_{\mathfrak{P}}(\chi_{\alpha})$) if, and only if, $\{\chi_{\alpha}\}_{\alpha\in A}$ is complete in $H$ and that, when $\{\chi_{\alpha}\}_{\alpha\in A}$ is semicomplete in $H$ or in $L^{2}(G),$ both $H_{\mathfrak{P}}(\chi_{\alpha})$ and $H_{\mathfrak{P}}'(\chi_{\alpha})$ are non-zero: an example may be seen from using the \textit{standard Riemann-Lebesgue orthonormal set.} In general, we have the following.

\textbf{3.4 Lemma.} \textit{Let $(H,\langle\cdot,\cdot\rangle)$ denote any Hilbert space. Then $$H(\chi_{\alpha})\subseteq H_{\mathfrak{P}}'(\chi_{\alpha})$$ for any semicomplete orthonormal set $\{\chi_{\alpha}\}_{\alpha\in A}$ in $H.$}

\textbf{Proof.} Choose any $x\in H(\chi_\alpha),$ then $x=\sum_{\alpha\in A}\langle x,\chi_{\alpha}\rangle\chi_{\alpha}.$ Now if $\langle x,\chi_{\alpha}\rangle=0,$ for every $\alpha\in A,$ then $$x=\sum_{\alpha\in A}\langle x,\chi_{\alpha}\rangle\chi_{\alpha}=\sum_{\alpha\in A}(0)\chi_{\alpha}=0;$$ showing that $x=0$ as required$.\;\Box$

We shall refer to $H_{\mathfrak{P}}'(\chi_{\alpha})$ as the \textit{prime-Parseval subspace} of $H$ and the choice of this term is further reinforced by the following facts.

\textbf{3.5 Lemma.} (\textit{cf.} Lemma 2.2) \textit{Let $\{\chi_{\alpha}\}_{\alpha\in A}$ denote a semicomplete orthonormal set in a Hilbert space $(H,\langle\cdot,\cdot\rangle)$ and let $x\in H.$ Then $x\in H_{\mathfrak{P}}'(\chi_{\alpha})$ whenever $\parallel x\parallel^{2}=\sum_{\alpha\in A}\mid\langle x,\chi_{\alpha}\rangle\mid^{2}.$}

\textbf{Proof.} If $\parallel x\parallel^{2}=\sum_{\alpha\in A}\mid\langle x,\chi_{\alpha}\rangle\mid^{2}$ and $\mid\langle x,\chi_{\alpha}\rangle\mid=0$ (for every $\alpha \in A$), then $\parallel x\parallel^{2}=\sum_{\alpha\in A}\mid\langle x,\chi_{\alpha}\rangle\mid^{2}=\sum_{\alpha\in A}(0)=0;$ showing that $x=0.$ Hence $x\in H_{\mathfrak{P}}'(\chi_{\alpha}).\;\Box$

Lemma $3.5$ shows the first partial connection between the satisfaction of \textit{Parseval equality,} on one hand, and membership in the \textit{prime-Parseval subspace,} on the other. The last Lemma may also be seen as saying that the subset of $H$ given as $$\{x\in H:\;\parallel x\parallel^{2}=\sum_{\alpha\in A}\mid\langle x,\chi_{\alpha}\rangle\mid^{2},\;\mbox{for any orthonormal set $\{\chi_{\alpha}\}_{\alpha\in A}$}\}$$ is also a subset of $H_{\mathfrak{P}}'(\chi_{\alpha}),$ with clear equality when $\{\chi_{\alpha}\}_{\alpha\in A}$ is complete. It will be satisfying to also have the reverse inclusion, $$H_{\mathfrak{P}}'(\chi_{\alpha})\subseteq\{x\in H:\;\parallel x\parallel^{2}=\sum_{\alpha\in A}\mid\langle x,\chi_{\alpha}\rangle\mid^{2},\;\mbox{for any orthonormal set $\{\chi_{\alpha}\}_{\alpha\in A}$}\}$$ due to the importance of the Parseval equality in the fine properties of Fourier transform. We shall deal with this concern in Lemma $3.12.$

Even though a semicomplete orthonormal set $\{\chi_{\alpha}\}_{\alpha\in A}$ in $L^{2}(G)$ (or in a Hilbert space $(H,\langle,\cdot,\rangle)$) may not be dense, as it is generally expected of a complete orthonormal set, we may still however employ this orthonormal set to construct some dense subspaces of $L^{2}(G)$ (or of a Hilbert space $(H,\langle,\cdot,\rangle)$) as follows. Indeed, the following results on the \textit{Fourier subspace} for $L^{2}(G)$ are also valid for an arbitrary Hilbert space, $(H,\langle,\cdot,\rangle)$ and for a \textit{relative semicomplete orthonormal set} in $H.$

\textbf{3.6 Theorem.} \textit{Let $G$ denote a compact and let $\{\chi_{\alpha_{i}^{j}}\}_{{\alpha_{i}^{j}}\in A}$ denote a semicomplete orthonormal set on $G.$ Then $L^{2}(G)(\chi_{{\alpha_{i}^{j}}})$ is topologically dense in $L^{2}(G).$}

\textbf{Proof.} Since every $f\in L^{2}(G)$ may be expanded as $$\sum_{\lambda\in\widehat{G}}d(\lambda)\sum_{i,j=1}^{d(\lambda)}\langle f,u_{ij}^{\lambda}\rangle u_{ij}^{\lambda}$$ (with convergence in the norm of $L^{2}(G)$) it follows that for $\epsilon>0$ we have $$\parallel f-\sum_{\lambda\in\widehat{G}}d(\lambda)\sum_{i,j=1}^{d(\lambda)}\langle f,u_{ij}^{\lambda}\rangle u_{ij}^{\lambda}\parallel_{2}<\frac{\epsilon}{2}.$$ Now $$\parallel f-\sum_{j=1}^{n}\gamma_{j}\sum_{i=1}^{n}\beta_{ij}\sum_{\alpha_{i}^{j}\in A}\langle f,\chi_{\alpha_{i}^{j}}\rangle\chi_{\alpha_{i}^{j}}\parallel_{2}\leq\;\parallel f-\sum_{\lambda\in\widehat{G}}d(\lambda)\sum_{i,j=1}^{d(\lambda)}\langle f,u_{ij}^{\lambda}\rangle u_{ij}^{\lambda}\parallel_{2}$$ $$+\parallel \sum_{\lambda\in\widehat{G}}d(\lambda)\sum_{i,j=1}^{d(\lambda)}\langle f,u_{ij}^{\lambda}\rangle u_{ij}^{\lambda}-\sum_{j=1}^{n}\gamma_{j}\sum_{i=1}^{n}\beta_{ij}\sum_{\alpha_{i}^{j}\in A}\langle f,\chi_{\alpha_{i}^{j}}\rangle\chi_{\alpha_{i}^{j}}\parallel_{2}<\frac{\epsilon}{2}+\frac{\epsilon}{2}=
\epsilon.\;\Box$$

In more specific terms we have the following.

\textbf{3.7 Corollary.} \textit{Let $G$ denote a compact group and let $\{\chi_{\alpha_{i}^{j}}\}_{{\alpha_{i}^{j}}\in A}$ denote a semicomplete orthonormal set on $G.$ Then every $f\in L^{2}(G)$ can be expanded as $$f=\sum_{j=1}^{n}\gamma_{j}\sum_{i=1}^{n}\beta_{ij}\sum_{\alpha_{i}^{j}\in A}\langle f,\chi_{\alpha_{i}^{j}}\rangle\chi_{\alpha_{i}^{j}}$$ for some $\gamma_{j},\beta_{ij}\in\mathbb{C}$ with convergence in the norm on $L^{2}(G).\;\Box$}

We may refer to the expansion of $f$ in Corollary $3.7$ as a \textit{semi-Fourier series expansion} for $f\in L^{2}(G)$ or $H$ with respect to $\{\chi_{\alpha_{i}^{j}}\}_{{\alpha_{i}^{j}}\in A}.$ A stronger form of Theorem $3.6$ carved in the form of the equivalence of Lemma $2.2$ and which generalizes the fact that a mutually orthonormal family $\{\chi_{\alpha}\}_{\alpha\in A}$ is complete (in a Hilbert space $(H,\langle\cdot,\cdot\rangle)$) if, and only if, $H(\chi_{\alpha})=H$ (\textit{cf.} Lemma $2.2$) is also possible when the mutually orthonormal family $\{\chi_{\alpha}\}_{\alpha\in A}$ is semicomplete in $H.$ We prove this below in the special case of $H=L^{2}(G).$

\textbf{3.8 Theorem.} \textit{Let $G$ denote a compact group and let $\{\chi_{\alpha_{i}^{j}}\}_{{\alpha_{i}^{j}}\in A}$ denote a mutually orthonormal set on $G$ whose Fourier subspace is denoted as $L^{2}(G)(\chi_{\alpha_{i}^{j}}).$ Then $L^{2}(G)(\chi_{\alpha_{i}^{j}})$ is topologically dense in $L^{2}(G)$ if, and only if, $\{\chi_{\alpha_{i}^{j}}\}_{{\alpha_{i}^{j}}\in A}$ is semicomplete.}

\textbf{Proof.} That $L^{2}(G)(\chi_{\alpha_{i}^{j}})$ is topologically dense in $L^{2}(G)$ if $\{\chi_{\alpha_{i}^{j}}\}_{{\alpha_{i}^{j}}\in A}$ is semicomplete is the content of Theorem $3.6.$ Now choose $f\in L^{2}(G),$ then $$\parallel \sum_{\lambda\in\widehat{G}}d(\lambda)\sum_{i,j=1}^{d(\lambda)}\langle f,u_{ij}^{\lambda}\rangle u_{ij}^{\lambda}-\sum_{j=1}^{n}\gamma_{j}\sum_{i=1}^{n}\beta_{ij}\sum_{\alpha_{i}^{j}\in A}\langle f,\chi_{\alpha_{i}^{j}}\rangle\chi_{\alpha_{i}^{j}}\parallel_{2}$$ $$\leq\parallel \sum_{\lambda\in\widehat{G}}d(\lambda)\sum_{i,j=1}^{d(\lambda)}\langle f,u_{ij}^{\lambda}\rangle u_{ij}^{\lambda}-f\parallel_{2}\;+\;\parallel f-\sum_{j=1}^{n}\gamma_{j}\sum_{i=1}^{n}\beta_{ij}\sum_{\alpha_{i}^{j}\in A}\langle f,\chi_{\alpha_{i}^{j}}\rangle\chi_{\alpha_{i}^{j}}\parallel_{2}$$ $\frac{\epsilon}{2}+\frac{\epsilon}{2}=\epsilon\mbox{(using the \textit{Peter-Weyl theorem} and Corollary $3.7,$ respectively)}.\;\Box$

This Theorem would enable us to see the \textit{Peter-Weyl series expansion} of every $f\in L^{2}(G),$ given as $$f=\sum_{\lambda\in\widehat{G}}d(\lambda)\sum_{i,j=1}^{d(\lambda)}\langle f,u_{ij}^{\lambda}\rangle u_{ij}^{\lambda}$$ (with convergence in the $L^{2}-$norm), as the restriction of the \textit{semi-Fourier series expansion} $$f=\sum_{j=1}^{n}\gamma_{j}\sum_{i=1}^{n}\beta_{ij}\sum_{\alpha_{i}^{j}\in A}\langle f,\chi_{\alpha_{i}^{j}}\rangle\chi_{\alpha_{i}^{j}}$$ to the \textit{standard Peter-Weyl (complete) mutually orthonormal set} $\{\sqrt{d(\lambda)}u_{ij}^{\lambda}\}.$ Indeed Theorem $3.8$ leads to the same conclusion for the \textit{prime-Parseval subspace} $L^{2}(G)_{\mathfrak{P}}'(\chi_{\alpha_{i}^{j}}).$

\textbf{3.9 Corollary.} \textit{Let $G$ denote a compact group and let $\{\chi_{\alpha_{i}^{j}}\}_{{\alpha_{i}^{j}}\in A}$ denote a mutually orthonormal set on $G.$ Then $L^{2}(G)_{\mathfrak{P}}'(\chi_{\alpha_{i}^{j}})$ is topologically dense in $L^{2}(G)$ if, and only if, $\{\chi_{\alpha_{i}^{j}}\}_{{\alpha_{i}^{j}}\in A}$ is semicomplete.}

\textbf{Proof.} Consider Lemma $3.4$ in the light of Theorem $3.8.\;\Box$

The inclusion $L^{2}(G)(\chi_{\alpha_{i}^{j}})\subseteq L^{2}(G)_{\mathfrak{P}}'(\chi_{\alpha_{i}^{j}})$ of Lemma $3.4,$ when combined with both Theorem $3.7$ and Corollary $3.9,$ implies the following.

\textbf{3.10 Corollary.} \textit{$L^{2}(G)(\chi_{\alpha_{i}^{j}})$ is topologically dense in $L^{2}(G)_{\mathfrak{P}}'(\chi_{\alpha_{i}^{j}}).\;\Box$}

The converse of Lemma $3.5$ is now immediate for both $L^{2}(G)_{\mathfrak{P}}'(\chi_{\alpha_{i}^{j}})$ and (even) $H_{\mathfrak{P}}'(\chi_{\alpha_{i}^{j}})$ in any arbitrary Hilbert space, $(H,\langle\cdot,\cdot\rangle).$

\textbf{3.11 Lemma.}(\textit{cf.} Lemma $2.2$) \textit{Let $G$ denote a compact group and let $\{\chi_{\alpha}\}_{{\alpha}\in A}$ denote a mutually orthonormal set on $G.$ Then $f\in L^{2}(G)_{\mathfrak{P}}'(\chi_{\alpha})$ if, and only if, $\parallel f\parallel_{2}^{2}=\sum_{\alpha\in A}\mid\langle f,\chi_{\alpha}\rangle\mid^{2}.$}

\textbf{Proof.} Let $f\in L^{2}(G)_{\mathfrak{P}}'(\chi_{\alpha}).$ We may take $f\in L^{2}(G)(\chi_{\alpha})$ due to Corollary $3.10;$ so that $f=\sum_{\alpha\in A}\langle f,\chi_{\alpha}\rangle\chi_{\alpha}.$ Hence $$0=\parallel f-\sum_{\alpha\in A}\langle f,\chi_{\alpha}\rangle\chi_{\alpha}\parallel_{2}^{2}=\parallel f\parallel_{2}^{2}-\sum_{\alpha\in A}\mid \langle f,\chi_{\alpha}\rangle\mid^{2},$$ as required$.\;\Box$

Hence, the \textit{prime-Parseval subspace} $L^{2}(G)_{\mathfrak{P}}'(\chi_{\alpha_{i}^{j}})$ may finally be seen (for some orthonormal set $\{\chi_{\alpha_{i}^{j}}\}_{{\alpha_{i}^{j}}\in A}$) as $$L^{2}(G)_{\mathfrak{P}}'(\chi_{\alpha_{i}^{j}})=\{f\in L^{2}(G):\;\parallel f\parallel_{2}^{2}=\sum_{\alpha_{i}^{j}\in A}\mid\langle f,\chi_{\alpha_{i}^{j}}\rangle\mid^{2}\}$$

We now have enough preparation to introduce a Fourier transform $f\mapsto \widehat{f}$ on the \textit{prime-Parseval subspace,} $L^{2}(G)_{\mathfrak{P}}'(\chi_{\alpha_{i}^{j}}).$

Consider $f\in L^{2}(G)$ and for every $\alpha\in A$ define the matrix $\widehat{f}(\alpha)$ whose entries are given as $$\widehat{f}(\alpha)_{ij}:=\widehat{f}(\alpha_{i}^{j}).$$ That is, $\widehat{f}(\alpha)_{ij}:=\langle f,\chi_{\alpha_{i}^{j}}\rangle,$ for $1\leq i,j\leq n.$ The Parseval inequality of $L^{2}(G)_{\mathfrak{P}}'(\chi_{\alpha_{i}^{j}})$ (in Lemma $3.11$) therefore becomes $\parallel f\parallel_{2}^{2}=\sum_{\alpha\in A}\parallel\widehat{f}(\alpha)\parallel^{2},$ for every $f\in L^{2}(G)_{\mathfrak{P}}'(\chi_{\alpha_{i}^{j}}),$ where $\parallel\widehat{f}(\alpha)\parallel^{2}$ is the Hilbert-Schimdt norm of the matrix $$\widehat{f}(\alpha)=(\widehat{f}(\alpha)_{ij})_{i,j=1}^{n}=(\widehat{f}(\alpha_{i}^{j}))_{i,j=1}^{n}.$$
In other words, and in terms of our choice of indexing $A,$ we have $$\parallel f\parallel_{2}^{2}=\sum_{i,j=1}^{n}\sum_{\alpha_{i}^{j}\in A}\parallel\widehat{f}(\alpha_{i}^{j})\parallel^{2},$$ for $f\in L^{2}(G)_{\mathfrak{P}}'(\chi_{\alpha_{i}^{j}}).$

\textbf{3.12 Definition.} Set $L^{2}(A)$ as the space of matrix-valued functions $\varphi$ on $A$ with values in $\bigcup_{n=1}^{\infty}M_{n}(\mathbb{C})$ satisfying\\
$(i)$ $\varphi(\alpha_{i}^{j})\in M_{n}(\mathbb{C})$ for all $\alpha_{i}^{j}\in A$ and\\
$(ii)$ $\sum_{i,j=1}^{n}\sum_{\alpha_{i}^{j}\in A}\parallel\varphi(\alpha_{i}^{j})\parallel^{2}<\infty.\;\Box$

The inner product $(\cdot,\cdot)$ on $L^{2}(A)$ given as $$(\varphi,\psi):=\sum_{i,j=1}^{n}\sum_{\alpha_{i}^{j}\in A}tr(\varphi(\alpha_{i}^{j})\psi(\alpha_{i}^{j})^{*}),$$ $\varphi,\psi\in L^{2}(A)$ converts $(L^{2}(A),(\cdot,\cdot))$ into a Hilbert space. We can then establish a connection between the \textit{prime-Parseval subspace} $L^{2}(G)_{\mathfrak{P}}'(\chi_{\alpha_{i}^{j}})$ (which is a Hilbert subspace of $L^{2}(G)$) and $L^{2}(A).$

\textbf{3.13 Theorem.} (Fourier image of the \textit{prime-Parseval subspace}) \textit{Let $G$ denote a compact group and let $\{\chi_{\alpha_{i}^{j}}\}_{{\alpha_{i}^{j}}\in A}$ denote a semicomplete mutually orthonormal set on $G.$ Then the map $$\mathcal{H}:L^{2}(G)_{\mathfrak{P}}'(\chi_{\alpha_{i}^{j}})\rightarrow L^{2}(A):f\mapsto\mathcal{H}(f):=\widehat{f}$$ is an isometry of $L^{2}(G)_{\mathfrak{P}}'(\chi_{\alpha_{i}^{j}})$ onto $L^{2}(A).\;\Box$}

Theorem $3.13$ is very familiar when the semicomplete mutually orthonormal set $\{\chi_{\alpha_{i}^{j}}\}_{{\alpha_{i}^{j}}\in A}$ is the complete mutually orthonormal set $\{\sqrt{d(\lambda)}u_{ij}^{\lambda}\}.$ We do not yet know the general connection between the set $A$ and the dual group $\widehat{G},$ except in the special cases of the \textit{standard Riemann-Lebesgue (semicomplete) orthonormal sets.} We however see $A$ as a general form of $\widehat{G}$ which may take the usual form of $\widehat{G}$ in specific cases. If we set $$H_{i}^{\alpha}:=\sum_{j=1}^{n}\mathbb{C}\chi_{\alpha_{i}^{j}},$$ for $\alpha=\alpha_{i}^{j}\in A$ and $i\in\{1,\cdots,n\},$ then the Hilbert subspace $L^{2}(G)_{\mathfrak{P}}'(\chi_{\alpha_{i}^{j}})$ of $L^{2}(G)$ has the direct-sum decomposition $$L^{2}(G)_{\mathfrak{P}}'(\chi_{\alpha})=\bigoplus_{\alpha\in A}\bigoplus_{i=1}^{n}H_{i}^{\alpha}.$$

The results of this section laid a foundation for harmonic analysis of the \textit{prime-Parseval subspace} $H_{\mathfrak{P}}'(\chi_{\alpha_{i}^{j}})$ with respect to a semicomplete orthonormal set $\{\chi_{\alpha_{i}^{j}}\}_{{{\alpha_{i}^{j}}}\in A}$ in a Hilbert space, $H.$ Having considered the case of the Hilbert space $L^{2}(G),$ for a compact group $G,$ in this section it will a delight to use these foundational results (on both $H_{\mathfrak{P}}'(\chi_{\alpha_{i}^{j}})$ and $L^{2}(G)_{\mathfrak{P}}'(\chi_{\alpha_{i}^{j}})$) in the understanding of further properties of $L^{2}(G)_{\mathfrak{P}}'(\chi_{\alpha_{i}^{j}})$ in the full sight of the semicompleteness of $\{\chi_{\alpha_{i}^{j}}\}_{{{\alpha_{i}^{j}}}\in A}.$ We shall give a very short introduction to this type of study for a connected semisimple Lie group in the next section.

It is clear from Lemma $3.2,$ for \textit{standard (Riemann-Lebesgue)} examples of a semicomplete orthonormal  set in an arbitrary Hilbert space $(H,\langle\cdot,\cdot\rangle)$ or in $L^{2}(G),$ that the non-zero constants $\gamma_{j}$ and $\beta_{ij}$ would always be $\gamma_{j}=\beta_{ij}=1$ for $1\leq i,j\leq\mid\widehat{G}\setminus\{\lambda_{0}^{(1)},\lambda_{0}^{(2)},\cdots\}\mid.$ However, for \textit{non-standard} examples of a semicomplete orthonormal  set in an arbitrary Hilbert space $(H,\langle\cdot,\cdot\rangle)$ or in $L^{2}(G),$ the \textit{semi-Fourier series expansion} of Corollary $3,7$ may have to be broken down in order for general expressions for $\gamma_{j}$ and $\beta_{ij}$ to be known. A first result along this line is the following.

\textbf{3.14 Lemma.} \textit{Let $\{\chi_{\alpha_{i}^{j}}\}_{\alpha_{i}^{j}\in A}$ denote a semicomplete orthonormal set in a Hilbert space $(H,\langle\cdot,\cdot\rangle)$ and let $x\in H.$ Then $$\langle x,\chi_{\alpha_{i}^{i}}\rangle=\gamma_{i}\beta_{ii}\langle x,\chi_{\alpha_{i}^{i}}\rangle,$$ for  $1\leq i\leq n.$ In particular, $\gamma_{i}\beta_{ii}=1.$}

\textbf{Proof.} We have that $\langle x,\chi_{\alpha_{k}^{l}}\rangle=\sum_{j=1}^{n}\gamma_{j}\sum_{i=1}^{n}\beta_{ij}\sum_{\alpha_{i}^{j}\in A}\langle x,\chi_{\alpha_{i}^{j}}\rangle\langle\chi_{\alpha_{i}^{j}},\chi_{\alpha_{k}^{l}}\rangle.$ Due to the orthogonality of the set $\{\chi_{\alpha_{i}^{j}}\}_{\alpha_{i}^{j}\in A}$ the above equality reduces to $\langle x,\chi_{\alpha_{i}^{i}}\rangle=\gamma_{i}\beta_{ii}\langle x,\chi_{\alpha_{i}^{i}}\rangle,$ for  $1\leq i\leq n$ as required.

Now $(1-\gamma_{i}\beta_{ii})\langle x,\chi_{\alpha_{i}^{i}}\rangle=0$ from where we have $\gamma_{i}\beta_{ii}=1.\;\Box$

\ \\
\ \\
{\bf \S 4. K-semicomplete orthonormal set in a semisimple Lie group.}

The success in $\S 3.$ of the use of the notion of a \textit{semicomplete orthonormal} set in the harmonic analysis of a compact group, culminating in the extraction and elucidation of the \textit{prime-Parseval subsapce} as well as its Fourier image, shows the central importance and the correct us of \textit{Parseval equality} and the concept of \textit{completeness} (of an orthonormal set) in the abstract Peter-Weyl theory of a compact group and in the understanding of the hitherto unknown subspaces of $L^{2}(G)$ under the influence of the Fourier transform. This study (which led us to the consideration of the \textit{prime-Parseval subspace} $L^{2}(G)_{\mathfrak{P}}'(\chi_{\alpha_{i}^{j}})$ corresponding to a semicomplete orthonormal set $\{\chi_{\alpha_{i}^{j}}\}_{{{\alpha_{i}^{j}}}\in A}$ on $G$) is reminiscence of and may be compared with the extraction and harmonic analysis of the \textit{Schwartz algebra} in the $L^{2}-$theory of semisimple Lie groups which was started in the Yale thesis $[1(a.)]$ of James Arthur (continued and completed in two later manuscripts, $[1(b.)]$ and $[1(c.)]$). In a more recent publication, harmonic analysis of other spaces of functions on semisimple Lie groups, namely of the space of \textit{spherical convolutions,} has been introduced in $[3.]$ leading to the explicit construction of the corresponding \textit{Plancherel formula} for such functions. The present paper has also introduced the \textit{Fourier} and \textit{prime-Parseval subspaces} of $L^{2}(G)$ (or of any arbitrary Hilbert space, $(H,\langle\cdot,\cdot\rangle)$).

Having shown in $\S 3.$ the essential importance of the Parseval equality (which is the precursor of the Plancherel formula) in the consideration of the actual subspace of $L^{2}(G)$ under the natural action of the Fourier transform, we shall here consider studying the same theory (of a semicomplete orthonormal set) but for all semisimple Lie groups, having removed the impediments posed by the \textit{completeness} for orthonormal sets on such Lie groups.

It is well-known that orthonormal sets (of functions and polynomials) are numerous and readily available in the $L^{2}-$space (and more recently in some distinguished subspaces of the $L^{2n}-$spaces $[4.]$) of semisimple Lie groups. Indeed every semisimple Lie group has its corresponding orthonormal set, an example is $G=SL(2,\R)$ and its \textit{Legendre functions.}

Even though these sets of orthonormal functions and polynomials are central to harmonic analysis on these groups, their direct importance in or contribution to the decomposition of (sub-)spaces of $L^{2}(G)$ or expansion of their members is not yet known. In the outlook of the present section (and of the entire paper) any orthonormal set on a semisimple Lie group known to have been \textit{$K-$semicomplete} (in the sense to be soon made precise) could be a basis of some subspaces of $L^{2}(G).$

\textbf{4.1 Definition.} ($K-$semicomplete orthonormal set) \textit{Let $G=KAN$ denote the Iwasawa decomposition of a connected semisimple Lie group $G$ with finite center. An orthonormal set $\{\chi_{\alpha}\}_{{{\alpha}\in A}}$ on $G$ is said to be $K-$semicomplete whenever its restriction to $K,$ written as $\{(\chi_{\alpha})_{\mid_{K}}\}_{{{\alpha}}\in A},$ is a semicomplete orthonormal set in $L^{2}(K).\;\Box$}

It is relatively easy to construct a \textit{$K-$semicomplete orthonormal set} on any connected semisimple Lie group $G,$ from any given semicomplete orthonormal set on $K$ as follows.

\textbf{4.2 An example.} Choose any of the numerous orthonormal sets $\{\xi_{\alpha}\}_{\alpha\in A}$ in $L^{2}(K)$ as constructed in $\S 3.$ and, for every $x=kan\in G,$ define the map $\chi_{\alpha}:G\rightarrow\C$ as $$\chi_{\alpha}(x)=\chi_{\alpha}(kan):=e^{f(an)}\xi_{\alpha}(k),$$
where $f:AN\rightarrow\C$ satisfies\\
$(i)$ $f(1)=0,$\\
$(ii)$ $\int_{AN}e^{2\Re(f(an))}dadn=1$ and\\
$(iii)$ $\int_{AN}g(kan)(e^{\overline{f(an)}+f(a_{1}n_{1})})dadn=g(k),$ for $g\in L^{2}(G),$ $a_{1}\in A,$ $n_{1}\in N$ and the normalized Haar measures $da$ and $dn$ on $A$ and $N,$ respectively.

\textbf{Proof.} Observe that since $$\chi_{\alpha}(x)=\chi_{\alpha}(kan):=e^{f(an)}\xi_{\alpha}(k),$$ then for any $k\in K$ $$\chi_{\alpha}(k)=\chi_{\alpha}(k\cdot1\cdot1):=e^{f(1\cdot1)}\xi_{\alpha}(k)=\xi_{\alpha}(k).$$ For any $\alpha_{1},\alpha_{2}\in A,$ we have $$\langle\chi_{\alpha_{1}},\chi_{\alpha_{2}}\rangle=
\int_{K}(\int_{AN}e^{2\Re(f(an))}dadn)\xi_{\alpha_{1}}(k)\overline{\xi_{\alpha_{2}}(k)}dk=
\langle\xi_{\alpha_{1}},\xi_{\alpha_{2}}\rangle$$ and $$\parallel\chi_{\alpha}\parallel_{2}^{2}=
\int_{K}(\int_{AN}e^{2\Re(f(an))}dadn)\mid\xi_{\alpha_{1}}(k)\mid^{2}dk=
\parallel\xi_{\alpha}\parallel_{2}^{2}=1;$$ showing that $\{\chi_{\alpha}\}_{\alpha\in A}$ is an orthonormal set on $G.$ Its $K-$semicompleteness is also shown as follows. For a pre-assigned $\epsilon>0,$ we have that $$\parallel \sum_{\lambda\in\widehat{G}}d(\lambda)\sum_{i,j=1}^{d(\lambda)}\langle g,u_{ij}^{\lambda}\rangle u_{ij}^{\lambda}-\sum_{j=1}^{n}\gamma_{j}\sum_{i=1}^{n}\beta_{ij}\sum_{\alpha_{i}^{j}\in A}\langle g,\chi_{\alpha_{i}^{j}}\rangle\chi_{\alpha_{i}^{j}}\parallel_{2}$$
$$=\parallel \sum_{\lambda\in\widehat{G}}d(\lambda)\sum_{i,j=1}^{d(\lambda)}\langle g,u_{ij}^{\lambda}\rangle u_{ij}^{\lambda}-\int_{K}[\int_{AN}g(kan)(e^{\overline{f(an)}+f(a_{1}n_{1})})dadn]\cdot$$ $$\sum_{j=1}^{n}\gamma_{j}\sum_{i=1}^{n}\beta_{ij}\sum_{\alpha_{i}^{j}\in A}\overline{\xi_{\alpha_{i}^{j}}(k)}dk\;\xi_{\alpha_{i}^{j}}\parallel_{2}$$
$$=\parallel \sum_{\lambda\in\widehat{G}}d(\lambda)\sum_{i,j=1}^{d(\lambda)}\langle g,u_{ij}^{\lambda}\rangle u_{ij}^{\lambda}-\sum_{j=1}^{n}\gamma_{j}\sum_{i=1}^{n}\beta_{ij}\sum_{\alpha_{i}^{j}\in A}\langle g,\xi_{\alpha_{i}^{j}}\rangle\xi_{\alpha_{i}^{j}}\parallel_{2}<\epsilon.\;\Box$$

For any $K-$semicomplete orthonormal set $\{\chi_{\alpha}\}_{\alpha\in A}$ on $G$ the corresponding \textit{Fourier subspace} $L^{2}(G)(\chi_{\alpha})$ of $L^{2}(G)$ is also given as $$L^{2}(G)(\chi_{\alpha}):=\{f\in L^{2}(G):\;f=\sum_{\alpha\in A}\langle f,\chi_{\alpha}\rangle \chi_{\alpha}\}$$ while the \textit{prime-Parseval subspace} is $$L^{2}(G)_{\mathfrak{P}}'(\chi_{\alpha}):=\{f\in L^{2}(G):\;\langle f,\chi_{\alpha}\rangle=0\;\mbox{(for every $\alpha\in A$) implies}\;f=0\}.$$

Clearly $L^{2}(K)_{\mathfrak{P}}'(\sqrt{d(\lambda)}u_{ij}^{\lambda})=L^{2}(K)$ (from Lemma $2.2\;(ii)$), both subspaces $L^{2}(K)(\chi_{\alpha})$ and $L^{2}(K)_{\mathfrak{P}}'(\chi_{\alpha})$ are topologically dense in $L^{2}(K)$ (from Theorems $3.6$ and $3.8$ and Corollary $3.9$) and there exists an isometry of $L^{2}(K)_{\mathfrak{P}}'(\chi_{\alpha})$ onto $L^{2}(A)$ (from Theorem $3.13$). We shall resume the study of the subspaces $L^{2}(G)(\chi_{\alpha})$ and $L^{2}(G)_{\mathfrak{P}}'(\chi_{\alpha})$ (for connected semisimple Lie groups, $G$) in another paper.
\ \\
\ \\
{\bf   References.}
\begin{description}
\item [{[1.]}] Arthur, J. G., $(a.)$ \textit{Harmonic analysis of tempered distributions on semisimple Lie groups of real rank one}, Ph.D. Dissertation, Yale University, $1970;$ $(b.)$ \textit{Harmonic analysis of the Schwartz space of a reductive Lie group I,} mimeographed note, Yale University, Mathematics Department, New Haven, Conn; $(c.)$ \textit{Harmonic analysis of the Schwartz space of a reductive Lie group II,} mimeographed note, Yale University, Mathematics Department, New Haven, Conn.
    \item [{[2.]}] Gangolli, R. and Varadarajan, V. S., \textit{Harmonic analysis of spherical functions on real reductive groups,} Ergebnisse der Mathematik und iher Genzgebiete, $vol.$ {\bf 101}, Springer-Verlag, Berlin-Heidelberg. $1988.$
                \item [{[3.]}] Oyadare, O. O., On harmonic analysis of spherical convolutions on semisimple Lie groups, \textit{Theoretical Mathematics and Applications}, $vol.$ $\textbf{5},$ no.: {\bf 3}. ($2015$), pp. 19-36.
                    \item [{[4.]}] Oyadare, O. O., Hilbert-substructure of real measurable spaces on reductive Groups, $I;$ Basic Theory, \textit{J. Generalized Lie Theory Appl.,} $vol.$ $\textbf{10},$ Issue {\bf 1}. ($2016$).
                    \item [{[5.]}] Sugiura, M., \textit{Unitary representations and harmonic analysis - an introduction} North-Holland Mathematical Library, $vol.$ {\bf 44}, Kodansha Scientific Books, Tokyo. $1990.$
                        \end{description}
\end{document}